\documentclass[12pt]{article}
\usepackage{latexsym}
\usepackage{amsbsy}
\usepackage{graphicx}
\usepackage{graphics}
\usepackage[initials, nobysame]{amsrefs}
\usepackage{enumerate}

\usepackage{todonotes}
\usepackage{comment}

\usepackage{tikz}
\usepackage{pgf}

\usepackage{pgfplots}

\DefineSimpleKey{bib}{primaryclass}{}
\DefineSimpleKey{bib}{archiveprefix}{}

\BibSpec{arXiv}{%
  +{}{\PrintAuthors}{author}
  +{,}{ \textit}{title}
  +{}{ \parenthesize}{date}
  +{,}{ arXiv }{eprint}
  +{,}{ primary class }{primaryclass}
}

\usepackage{amssymb}

\PassOptionsToPackage{hyphens}{url}
\usepackage{hyperref}

\definecolor{greyish}{gray}{0.7}
\definecolor{verylightgrey}{gray}{0.9}
\definecolor{greenish}{rgb}{0,0.5,0}
\definecolor{junglegreen}{rgb}{0.16,0.67,0.53}
\definecolor{lightblue}{rgb}{0.5,0.5,1}
\definecolor{yellowish}{rgb}{1,1,0}
\definecolor{ochre}{rgb}{0.80,0.47,0.13}
\definecolor{lightviolet}{rgb}{0.75,0.5,1}
\definecolor{indigo}{rgb}{0.25,0,1}
\definecolor{barbiepink}{rgb}{0.85,0.09,0.52}

\title{Dominic Welsh: his work and influence}
\author{
Graham Farr\\
Department of Data Science and AI, Faculty of IT, \\
Monash University, Clayton, VIC 3800, Australia;
\\\texttt{Graham.Farr@monash.edu}   
\and
Dillon Mayhew\\
School of Computing, University of Leeds, Leeds, LS2 9JT, UK;
\\\texttt{d.mayhew@leeds.ac.uk}
\and
James Oxley\\
Mathematics Department, Louisiana State University, \\Baton Rouge, LA 70803-4918, USA; \\\texttt{oxley@math.lsu.edu}
}

\date{28 June 2024}

\begin{document}

\maketitle

\begin{abstract}
We review the work of Dominic Welsh (1938--2023), tracing his remarkable influence through his theorems, expository writing, students, and interactions.  He was particularly adept at bringing different fields together and fostering the development of mathematics and mathematicians.  His contributions ranged widely across discrete mathematics over four main career phases: discrete probability, matroids and graphs, computational complexity, and Tutte-Whitney polynomials.  We give particular emphasis to his work in matroid theory and Tutte--Whitney polynomials.
\end{abstract}


\begin{center}
\includegraphics[height=8cm]{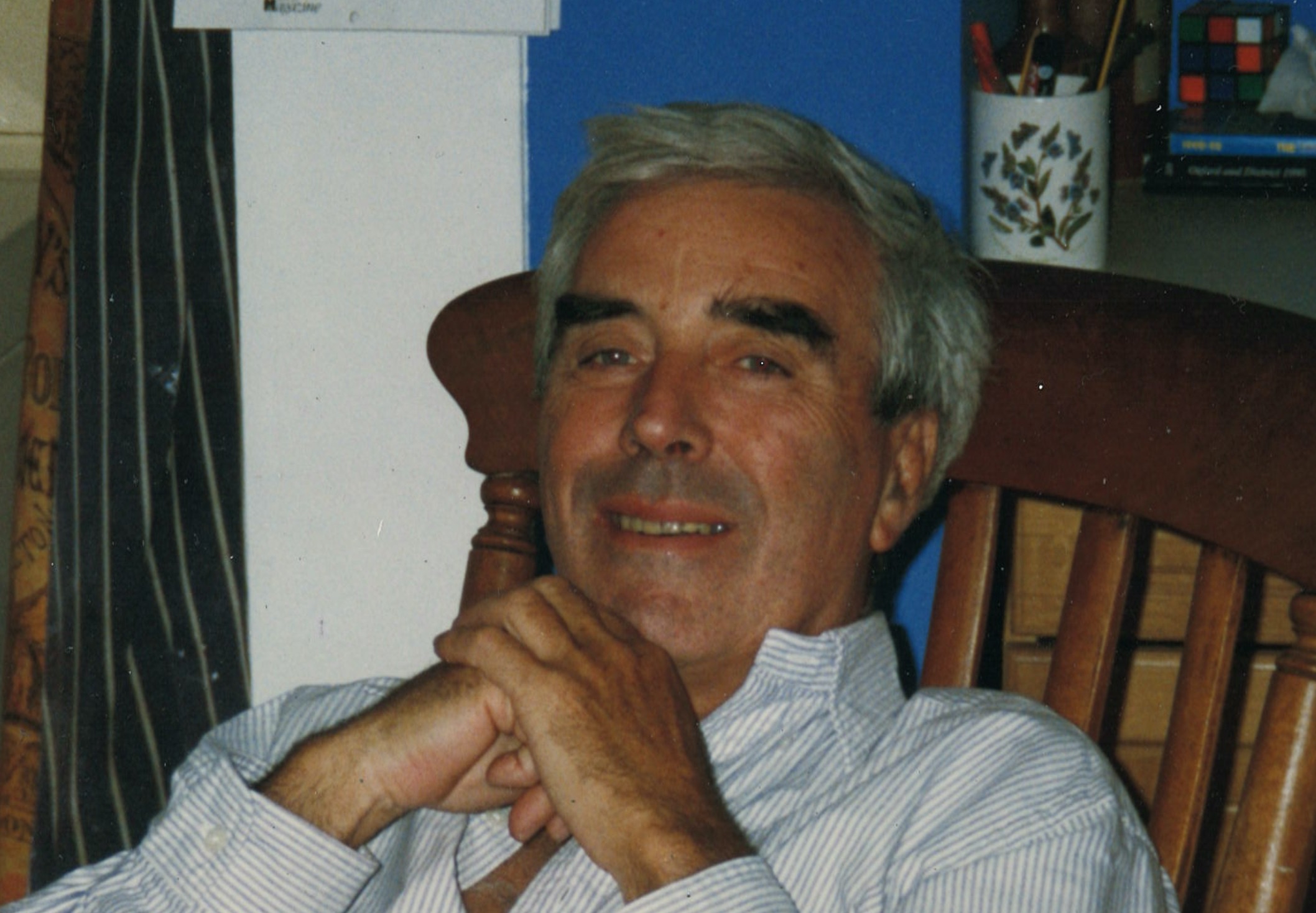}

\tiny Dominic Welsh at home in Oxford, mid- to late 1980s; photo by James Oxley

\end{center}

\section{Introduction}

The Discrete Mathematics community lost one of its greatest leaders and most treasured members with the death of Dominic Welsh on 30 November 2023.

Dominic had wide-ranging interests and many collaborators (57 listed in MathSciNet).  Most of his research was in at least one of three main themes:  discrete probability, matroids and graphs, and computational complexity.  These themes are captured well in the title \textit{Combinatorics, Complexity and Chance} of the tribute volume (edited by Geoffrey Grimmett and Colin McDiarmid) published by Oxford University Press in 2007, soon after he retired \cite{grimmett-mcdiarmid2007}.  That book is an excellent source of reflections on Dominic’s work.  It includes Oxley’s chapter on ``The contributions of Dominic Welsh to matroid theory'' \cite{oxley-07} which is a source of further information and reflections on that theme in Dominic’s work.  Here we will try to give an overview of his career,
indicating his range of interests and attempting to trace some of the sources of his remarkable influence.  We have divided his research career into four phases: probability, matroids, complexity, and Tutte polynomials.  The boundaries between these phases are not sharp, and each phase contained some seeds of later phases.

This article is an expanded version of our \textit{Matroid Union} post \cite{farr-mayhew-oxley2024} of February 2024.  Geoffrey Grimmett has written a wonderful obituary \cite{grimmett2024a} which includes a biography, reminiscences, and review of Dominic's contributions.  It also links to a list of Dominic's publications at \cite{grimmett2024b}.  Readers unfamiliar with Dominic's early work in discrete probability will find Grimmett's exposition of that area especially fascinating.  The present article is complementary through its areas of emphasis, its particular reflections on his influence, and its commentary on some work by Dominic's students.  (In discussing work by his students, we have had to be selective, focusing mostly on Dominic's main themes.)

\section{First phase: discrete probability}
\label{sec:discrete-prob}

Dominic's doctoral supervisor was John Hammersley (1920--2004), an eminent probability theorist and applied mathemtician.
Hammersley did not have a PhD/DPhil, though he was awarded Doctor of Science degrees later in his career by Cambridge and Oxford Universities.  Dominic recalled that Hammersley had a plain-speaking, direct manner in research supervision which worked well for him.  On one occasion, Dominic turned up to one of their meetings and confessed that he hadn’t done anything since their previous meeting.  Hammersley ended the meeting immediately and told him to come back when he’d done some work!  Geoffrey Grimmett writes of Hammersley that ``his unashamed love of a good problem has been an inspiration to many'' \cite{grimmett-89}*{p.~vii}, and the same can certainly be said of Dominic.

Dominic’s doctoral research was about percolation on square-lattice graphs \cite{welsh-64,hammersley-welsh1965}, in which information spreads out from the origin according to some local randomness on the edges or vertices.  
This topic had been pioneered in the late 1950s by Hammersley, initially with Simon Broadbent \cites{broadbent-hammersley1957,hammersley-57a, hammersley-57b}.  The original spark was a question by Broadbent on the porosity, for penetration by a gas, of carbon granules used in gas masks \citelist{\cite{hammersley-morton1954}*{pp.~68,74--75} \cite{hammersley-welsh1980}}.  In Broadbent and Hammersley's graph model \cite{broadbent-hammersley1957},
the edges allow information (or gas, or liquid,
or infection, \ldots) to pass, with some probability $p$, or not, with probability $1-p$.  The choices for
different edges are independent and identically
distributed.  The \textit{percolation probability} of the infinite square lattice is
the probability that information introduced at the origin can reach an infinite set of vertices.

In his DPhil thesis, Dominic introduced \textit{first-passage percolation}, which a graph theorist might think of as \textit{shortest path percolation} (although we must use infima, rather than minima, over sets of paths, since the lattice graphs are infinite).  Here, the time taken for information to pass along each edge is randomly distributed, and we study how quickly it moves from the origin to a distant vertex or to a line parallel to the $x$- or $y$-axis.  One of Dominic's theoretical contributions was to introduce the notion of a \textit{subadditive stochastic process}.  This is a family $(X_{s,t}:s,t\in\mathbb{N}\cup\{0\})$
of nonnegative real-valued random variables, such that $X_{s,t}$ has finite expectation and is
\begin{itemize}
    \item \textit{stationary}, that is, $X_{s,t}$ and $X_{s+k,t+k}$ are identically distributed for all $k\in\mathbb{N}$, and
    \item \textit{subadditive}, that is, $X_{r,t} \le X_{r,s} + X_{s,t}$
    whenever $r\le s\le t$.
\end{itemize}
We give an illustrative example from Dominic's thesis \citelist{\cite{welsh-64}*{Theorem 4.1.1} \cite{hammersley-welsh1965}}.
A random function $\omega$ assigns independent identically distributed weights to the edges of the lattice, with the weight of each edge representing
the time taken to go along it.  The length of a path is the sum of the weights of its edges.  Define $X_{s,t}=X_{s,t}(\omega)$ to be the infimum of the lengths of all paths from $(s,0)$ to $(t,0)$.
The family $(X_{s,t}:s,t\in\mathbb{N}\cup\{0\})$ is shown to be a subadditive process.

Dominic developed the theory of subadditive processes and proved the first subadditive ergodic theorem.  This enabled him to establish the existence of various fundamental limits for first-passage percolation.  A more detailed account of Dominic's work on percolation may be found in \cite{grimmett2024a}*{\S6}.

Dominic’s early work in probability included some graph theory, typically in the context of the lattice graphs on which percolation was studied \cite{bondy-welsh-66}.  He took up Hammersley’s interest in self-avoiding walks (paths starting at the origin with no repeated vertex) \cite{hammersley-welsh1962} and statistical-mechanical models on these graphs, and remained interested in discrete probability throughout his career.  His two research students who followed most closely in these footsteps were Geoffrey Grimmett (DPhil, 1974) and Peter Donnelly (DPhil, 1983), both of whom have since become Fellows of the Royal Society.  But, throughout his career, there were many other points of contact with his early interest in probability.  These included a 96-page survey \cite{welsh-70}, a textbook with Grimmett \cites{grimmett-welsh-86, grimmett-welsh-14}, papers with Oxley on generalisations of percolation \cites{oxley-welsh-79a, oxley-welsh-79b}, some classical percolation problems \cite{welsh-93a}, randomised algorithms \cites{welsh-83, petford-welsh-89, chavezlomeli-welsh1996} including Markov Chain Monte Carlo methods \cites{bartels-welsh-95, DVW-96}, and especially his work on the Tutte polynomial through its connections to statistical mechanics (via partition functions and the random-cluster model) \cites{welsh1990, welsh1997a, welsh-98, welsh1999, welsh-merino-00}.

\section{Second phase: matroid theory}
\label{sec:matroids}

In order to capture the notions of dependence common to linear algebra and graph theory, Whitney \cite{whitney1935} introduced matroids in 1935.  A \textit{matroid} $M$ consists of a finite set $E$ and a nonnegative integer-valued function $r$ on $2^E$ such that, for all subsets $X$ and $Y$ of $E$,
\begin{enumerate}[(i)]
\item $r(X)\le|X|$;
\item $r(X)\le r(Y)$ whenever $X\subseteq Y$;
\item $r$ is \textit{submodular}, that is, $r(X\cup Y)+r(X\cap Y) \le r(X) + r(Y)$.
\end{enumerate}
The set $E$ and the function $r$ are the \textit{ground set} and the \textit{rank function} of $M$.  Its \textit{independent sets} are those sets $I$ for which $r(I)=|I|$.  Matroids can be defined in numerous other ways including in terms of their independent sets, their \textit{bases} (maximal independent sets), and their \textit{circuits} (minimal dependent sets); see \cite{oxley-92,oxley-11,welsh-76}.  While Whitney introduced the term ``matroid'', equivalent structures were introduced contemporaneously by Nakasawa \cites{nakasawa1935,nakasawa1936a,nakasawa1936b}.

To obtain an example of a matroid, let $E$ be a finite set of vectors in a vector space and let $r(X)$ be the dimension of the subspace spanned by $X$.  As another example, let $E$ be the edge set of a graph $G$ and, for each subset $X$ of $E$, let $r(X)$ be the maximum number of edges in a forest in the induced subgraph $G[X]$.  This last matroid, $M(G)$, is called the \textit{cycle matroid} of $G$.

Dominic was first attracted to matroid theory by a seminar at Oxford by C.St.J.A. Nash-Williams in 1966.
The seminar was likely to have been similar to a conference talk Nash-Williams gave around that time \cite{NashWilliams1967} and showed how matroids could be used to give much easier proofs of some theorems in graph theory \cite{oxley-07}*{p.~235}.  This material is covered in \cite{welsh-76}*{\S\S8.3--8.4}.  Dominic later wrote that the seminar ``had a huge effect on my mathematical interests'' \cite{welsh2003}.

It was only ten years after that seminar that Dominic's classic book \textit{Matroid Theory} was published \cite{welsh-76}.  This period was one of intense activity and quickly established him as one of the worldwide leaders in the field.

The most obvious signs of this activity were Dominic’s publications on matroids.
In an early paper, he generalised to binary matroids the classical duality between Eulerian and bipartite planar graphs \cite{welsh-69a}.  He proved new lower bounds on the number of non-isomorphic matroids (\cite{welsh-69b}, and \cite{piff-welsh-71} with Piff).  Dominic attributed the rapid growth of interest in matroids from the mid-1960s to the discovery of transversal matroids (see his comments in \cite{welsh-76}*{p.~6}) and he made important contributions to that topic.  One of the most significant of these was a generalisation, using submodular functions, of theorems on transversals by Hall, Rado, Perfect, and Ore \cite{welsh-71a}. He also pinned down the essential role played by submodularity in that theorem.  He introduced matroids based on generalised transversals \cite{welsh-69c} and showed that every binary transversal matroid is graphic (with de Sousa \cite{desousa-welsh-72}).

Given his foundational work on subadditive stochastic processes in percolation theory, it is rather intriguing that submodular functions played such a key role in Dominic's early work on matroid theory.  After all, for real-valued functions on the subsets of a set, submodularity implies subadditivity.  The connection between Dominic's subadditivity and combinatorial submodularity is not close; in his first-passage percolation context, the subadditivity was of a more restricted kind and does not immediately extend to submodularity.  Nonetheless, his appreciation of subadditivity may have made him receptive to the power and centrality of submodularity in matroid theory.

Like the eminent Bill Tutte before him, when Dominic worked with matroids, he was heavily influenced by what had earlier been proved for graphs.  The stated purpose of his paper \emph{Matroids versus graphs} (with Harary \cite{harary-welsh-68}) was to help graph theorists ``to appreciate the important link between graph theory and matroids''.  This marks the beginning of Dominic's use of survey papers to both develop a field and to build connections to other areas of mathematics.

In July 1969, Dominic ran an influential combinatorics conference at Oxford that is now recognised as the first British Combinatorial Conference (BCC).  So, only three years after first meeting matroids, he was bringing combinatorics researchers together to advance the field.  In 1972, Dominic organised another combinatorics meeting in Oxford; it is now viewed as the 3rd BCC.

Dominic edited the proceedings of the 1969 Oxford conference \cite{welsh-71b}.  It contains his paper \emph{Combinatorial problems in matroid theory} \cite{welsh-71c}.  This provides a fascinating window into the state of discrete mathematics in the late 1960s, and clearly shows Dominic’s taste for attractive and challenging conjectures, as well as his long-lasting influence on mathematical research.  Several of the problems he proposed have been resolved, while many are still open.  Some have developed into important areas of research.  Dominic lists four problems concerning the enumeration of matroids, two of which have been settled.  The asymptotic behaviour of the number of binary matroids has been established by Wild \cite{wild-05}. Lemos \cite{lemos-04} and, independently, Crapo and Schmitt \cite{crapo-schmitt-05} solved another problem by showing that the number of non-isomorphic matroids on $m+n$ elements is at least the product of the numbers of non-isomorphic matroids on $m$ and on $n$ elements.  Other problems in enumeration remain open and seem very difficult, despite significant advances by Rudi Pendavingh, Jorn van der Pol, Remco van der Hofstad, and Nikhil Bansal \cites{BPP-14, BPP-15, pendavingh-vanderpol-15, vdHPvdP-22}.

The climax of Dominic’s first ten years in matroid theory was the publication of his book, \textit{Matroid Theory}, by Academic Press in 1976 \cite{welsh-76}.  This quickly became the standard reference for the field.  It covered a much wider set of topics, and reached a much wider audience, than previous books in the field.  Gian-Carlo Rota described it as ``beautifully written and exceptionally thorough'' \cite{rota-83}.  Its style is compact and concise, exemplifying his own direct advice to students who wrote too verbosely:  ``cut the chat!'', ``if in doubt, cut it out!''.  At the same time, it is accessible and inviting.  Its influence has been enormous.

Given the way Dominic’s interests developed later, it is worth noting his book’s importance in promoting the Whitney rank generating function and the Tutte polynomial.  Among the earliest surveys of that topic were Chapter 15 of Dominic's book and the treatment given in Norman Biggs’s book \textit{Algebraic Graph Theory} \cite{biggs1974}.  These works helped to broaden the interest in polynomial invariants for graphs and matroids.

When Dominic was approached by Oxford University Press in the late 1980s about updating his book, he declined and recommended James Oxley to the OUP editor.  Oxley’s book was published in 1992 \cite{oxley-92} (second edition published in 2011 \cite{oxley-11}).  Dominic’s book was reprinted, with some corrections, by Dover in 2010 and continues to be a valuable resource.  Dominic’s generosity of spirit led to him expressing concern that the Dover reprint of his 1976 book would hurt sales of Oxley’s book.

Two of Dominic’s enumerative problems concerning unimodal sequences lie at the heart of a new and rapidly expanding effort to use techniques from algebraic geometry in discrete mathematics.  A unimodal sequence can have only one local maximum (which might be a single peak or a plateau).  In his 1969 Oxford conference paper, Dominic looked at two sequences of numbers that can be derived from any matroid, namely the number of flats of rank $r$ and the number of independent sets of size $r$, where $r$ ranges from 0 to the rank of the matroid.  Harper and Rota had conjectured that the first sequence is unimodal; this remains open \cites{rota-harper-71, rota-71}.  Dominic conjectured that the same is true for the second sequence.

The characteristic polynomial of a matroid is a natural generalisation of the chromatic polynomial of a graph.  Brylawski \cite{brylawski-82} and, independently, Lenz \cite{lenz-13} showed that the coefficients of the characteristic polynomial of a matroid can be related to the sequence of numbers of independent sets of some matroid.  In Chapter 15 of Dominic’s book, he conjectured that the absolute values of these coefficients are log-concave, thereby strengthening earlier unimodal conjectures of Rota and of Heron \cites{rota-71, heron-72}.  (Andrew Heron’s 1972 DPhil was supervised by Dominic.)  In a tour de force of convex geometry and commutative algebra, Adiprasito, Huh, and Katz proved Dominic’s log-concavity conjecture \cite{AHK-18}.  Indeed, the proof of the Heron--Rota--Welsh Conjecture is prominently noted in the citation for June Huh’s 2022 Fields Medal \cite{fields}.

A major part of Dominic’s influence was through his supervision of research students.  His first doctoral student was Adrian Bondy (DPhil, 1969) with whom he wrote a paper on transversal matroids \cite{bondy-welsh-71} and who became a leading graph theorist.  Dominic recalled with self-effacing mirth that one of the first problems he gave Bondy was nothing less than the famous, notoriously difficult and still-open Reconstruction Conjecture!  After Bondy, then Heron, came DPhils by Michael Piff (1972), Joan de Sousa, Geoffrey Grimmett (both 1974), Colin McDiarmid, Frank Dunstan (both 1975), Lawrence Matthews (1977), James Oxley (1978), and Paul Walton (1982), all in matroid theory except for Grimmett.

Many of these students worked on a variety of topics within the field.  Heron’s paper on matroid polynomials in the 3rd BCC \cite{heron-72} included his celebrated unimodal conjecture.  This paper seems to have been the first appearance of matroid polynomials in the work of one of Dominic’s students.  McDiarmid made many contributions over his time as a student, with emphasis on transversal matroids, their duals, and links with Menger’s Theorem; some ten of his papers are cited in Dominic’s book.  Two of Dominic’s papers with Oxley were his first that focused on the Tutte polynomial and linked it to his early interest in percolation \cites{oxley-welsh-79a, oxley-welsh-79b}.  Results in the first of these papers included a type of `Recipe Theorem' that generalised a theorem of Brylawski and gave conditions under which an invariant is essentially an evaluation of a Tutte polynomial.  The second paper contained a characterisation of Tutte invariants at the very general level of arbitrary clutters (Sperner systems) using a generalisation of percolation probability.

Dominic also developed a close collaboration with Paul Seymour, whose DPhil (1975) at Oxford was supervised by the prominent matroid theorist Aubrey Ingleton.  Some of the work with Seymour established new links, in both directions, between Dominic’s interests in probability and combinatorics.  In \cite{seymour-welsh-75} they used the Fortuin--Kasteleyn--Ginibre (FKG) inequality from statistical mechanics to prove new results in combinatorics.
In \cite{seymour-welsh-78}
they developed clever geometrical arguments in a proof of a correlation-type inequality now called the Russo-Seymour-Welsh lemma, or RSW lemma.  This has been enormously useful in percolation and related fields, and played a key role in Kesten's celebrated proof \cite{kesten-80} that the critical probability for bond percolation on the square lattice is $\frac{1}{2}$,
as recounted in \cite[\S6]{grimmett2024a}.

Dominic continued to work in matroid theory throughout his career.  But he was a voracious learner who was always searching for new challenges and for ways to tie together areas that may have initially seemed disparate.

\section{Third phase: computational complexity}
\label{sec:complexity}

Around the early 1980s, the emphasis of Dominic’s research and supervision shifted from matroid theory to computational complexity.  While this may seem like a big change, his interest in computation went back a long way.  During a spell at Bell Labs in 1960--1961 (part of one of Hammersley’s visits there), he did some programming using punched cards as was normal in that era \cite{welsh2006}.  Later, while doing his doctorate at Oxford, he supplemented some of his theoretical results with computational simulations using Monte Carlo methods at the Oxford University Computing Laboratory on its Ferranti Mercury (its ``first big electronic computer'' \cite{sufrin2007}), as recounted in the last chapter of his thesis \cite{welsh-64}*{Ch.~9}.  In further work on the spread of infection through a two-dimensional lattice graph, he used Monte Carlo methods on an Elliott 803 computer, a type that was then common in British universities and was often programmed using Algol \cite{morgan-welsh1965}.  He discussed the link between the percolation problems he worked on and critical-path problems for Program Evaluation and Review Technique (PERT) networks, an important topic in project management and operations research, and wrote two short Letters to the Editor in the journal \textit{Operations Research} (1965--1966) discussing computational aspects \cite{welsh1965,welsh1966}.  In 1967, with Oxford colleague Martin Powell, he published an important sequential graph-colouring algorithm, giving an algorithmic refinement of Brooks's Theorem \cite{welsh-powell-67}.  This is known as the Welsh--Powell algorithm and has been widely studied and implemented.  In a paper at the 5th BCC in 1975, his serious interest in complexity is evident in discussing the Aanderaa--Rosenberg Conjecture (as it was then) about how much of a graph’s adjacency matrix needs to be read in order to determine if the graph has some hereditary property \cite{milner-welsh-75}.  With Gordon Robinson (MSc, 1979), he proved early results on the computational complexity of matroid properties in a framework where a matroid is specified by an oracle for one of five common axiomatisations, namely via independent sets, bases, circuits, rank, or closure \cite{robinson-welsh-80}.

\begin{figure}
\begin{center}
\begin{tikzpicture}
\node[text width=2cm,align=center,black] () at (6.3,13) {PSPACE};
\draw[line width=0.1cm,black] (3.5,7) ellipse (4cm and 6.1cm);
\node[text width=2.2cm,align=center,blue] () at (5,9.4) {$\hbox{P}^{\hbox{\scriptsize\#P}}=\hbox{P}^{\hbox{\scriptsize PP}}$};
\draw[line width=0.1cm,blue] (3.5,5) ellipse (3.3cm and 4cm);
\node[text width=2cm,align=center,olive] () at (2.8,4.6) {\scriptsize BQP};
\draw[line width=0.1cm,olive] (3.5,2.7) ellipse (1.6cm and 1.7cm);
\node[text width=2cm,align=center,greenish] () at (2.8,3.6) {\scriptsize BPP};
\draw[line width=0.1cm,greenish] (3.5,2.25) ellipse (1.2cm and 1.2cm);
\node[text width=2cm,align=center,junglegreen] () at (3.8,2.95) {\scriptsize RP};
\draw[line width=0.1cm,junglegreen,rotate around={-34:(3.5,1.5)}] (3.5,2) ellipse (0.7cm and 0.9cm);
\node[text width=2cm,align=center,violet] () at (5.2,5.8) {NP};
\draw[line width=0.1cm,violet,rotate around={-34:(3.5,1.5)}] (3.5,3.55) ellipse (1.35cm and 2.5cm);
\node[text width=1cm,align=center,violet] () at (3.9,2) {\small P};
\draw[line width=0.1cm,violet] (3.5,1.5) ellipse (0.37cm and 0.42cm);
\node[text width=1cm,align=center,black] () at (3.65,1.48) {\tiny L};
\draw[very thick,black] (3.5,1.26) ellipse (0.13cm and 0.13cm);
\end{tikzpicture}
\end{center}
\caption{Complexity classes}
\label{fig:complexity}
\end{figure}
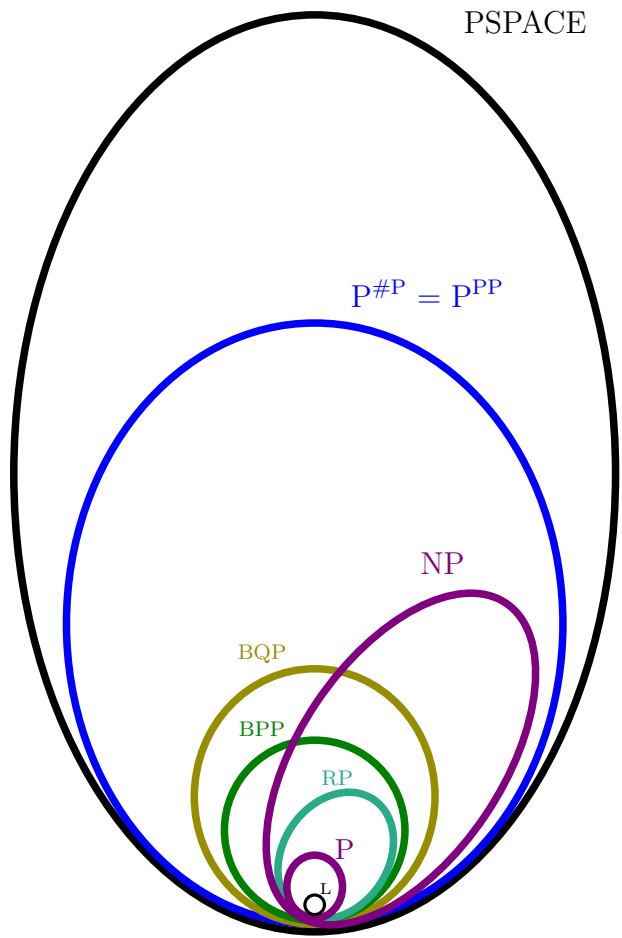

Dominic would often use Venn diagrams of complexity classes in his talks and conversations, and we do so now too.  Figure \ref{fig:complexity} shows some of the classes in which he was interested, presented much as he might have done.  The names of the classes are standard: see, for example, \citelist{\cite{garey-johnson1979} \cite{welsh-93b} \cite{talbot-welsh-10}  } for definitions and more information.  The classes LOGSPACE (abbreviated L) and PSPACE are included to frame the others, but they were not a major interest of Dominic's.  They are known to be distinct: LOGSPACE $\not=$ PSPACE, so the region PSPACE$\setminus$LOGSPACE in the diagram is nonempty.  But none of other regions indicated in this diagram is known to be nonempty.  In particular, none of the other subclass relations shown in the figure is known to be proper.  As well as the fundamental classes P and NP, Dominic was particularly interested in the probabilistic classes RP, BPP and PP, the quantum class BQP and the counting class \#P and its relatives.  The class $\hbox{NP}\,\cap\,\hbox{co-NP}$, and its relationship with other classes including P and RP, also intrigued him (though we have omitted it from Figure \ref{fig:complexity} for simplicity).

His DPhil students on computational complexity in the early to mid-1980s were Tony Mansfield (1982), Ken Regan (1986), Graham Farr (1986), and Keith Edwards (1986).  Tony Mansfield proved that determining the thickness of a graph is NP-hard \cite{mansfield1983}, solving a significant open problem in NP-completeness \cite[p.~286]{garey-johnson1979}.  Keith Edwards gave an impressive set of algorithms and hardness results for some fundamental existence and enumeration problems for colouring graphs of given minimum density or bounded genus \cite{edwards1986,edwards1992}.  As these examples illustrate, most of Dominic’s students in that era worked on specific graph-theoretic problems and classified their complexity (P versus NP-complete/hard etc.).  Ken Regan was the exception; he took the much more difficult path of proving new results about the complexity classes themselves and delving into the P-versus-NP problem.  He remains active in that field.  He and Richard Lipton have written about it in Lipton's blog, \textit{G\"odel’s Lost Letter and P=NP} \cite{lipton-regan}.

Dominic’s interest in computational complexity led naturally to an interest in cryptography, where complexity issues had come to the fore after the birth of public-key cryptography in the 1970s.  This led to a textbook, \textit{Codes and Cryptography} \cite{welsh-88} and supervision of some Masters students including Arun P. Mani (2004) and Douglas Stebila (2004).  Many years later, he would publish another textbook, with John Talbot, \textit{Complexity and Cryptography} \cite{talbot-welsh-10}.  Dominic’s attraction to complexity theory and cryptography was in spite of his recognition of the inherent difficulty of both subjects.

Dominic retained his earlier interests, supervising Peter Donnelly (1983) on discrete probability models on graphs and Manoel Lemos (1988) on matroid theory.  Some of his work with Donnelly  \cite{donnelly-welsh-83} inspired him to design a randomised 3-colouring algorithm which he studied with astrophysicist David Petford \cite{petford-welsh-89}.
That paper with Donnelly considered
antivoter models, in which the vertices of a graph are each given a colour Black or White, with the colour on a vertex randomly chosen to be biased away from whichever colour is more frequent on neighbouring vertices.  Some of Keith Edwards’s work also harked back to Dominic’s early interest in self-avoiding walks on square-lattice graphs \cite{edwards-85}.

This was an exceptionally busy period for Dominic, as he served variously as Chairman of the Board of the Faculty of Mathematics, Sub-Warden of Merton College, and Chairman of the British Combinatorial Committee, holding all three roles simultaneously for a time.  Somehow he also found the time to write two of his previously mentioned books.

\section{Fourth phase: Tutte--Whitney polynomials and their complexity}
\label{sec:tutte-polys}

The \textit{Whitney rank generating function} $R(M;x,y)$ of a matroid $M=(E,r)$ is defined by
\[
R(M;x,y)=\sum_{X\subseteq E} x^{r(E)-r(X)}y^{|X|-r(X)} .
\]
It was introduced by Whitney \cite{whitney1932} as a generalisation of the chromatic polynomial.  The \textit{Tutte polynomial} $T(M;x,y)$ is defined by
\[
T(M;x,y) = \left\{
\begin{array}{ll}
1,  &  \hbox{if $E=\emptyset$;}   \\
x\,T(M/e;x,y),  &  \hbox{if $e$ is a coloop;}   \\
y\,T(M\setminus e;x,y),  &  \hbox{if $e$ is a loop;}   \\
T(M\setminus e;x,y) + T(M/e;x,y),  &  \hbox{otherwise.}
\end{array}
\right.
\]
Here, coloops generalise bridges of graphs, and
the matroids $M\setminus e$ and $M/e$ are formed
using matroid generalisations of edge deletion and contraction in graphs.  These operations are commutative, so $T(M;x,y)$ is well defined.  The Tutte polynomial was introduced in \cite{tutte1947,tutte1948}.  The Whitney rank generating function and the Tutte polynomial are closely related \cite{tutte1947,tutte1954}:
\begin{equation}
\label{eq:tutte-from-whitney}
T(M;x,y) = R(M;x-1,y-1) .
\end{equation}

Evaluations of these polynomials at certain specific points or along certain curves in the $xy$-plane yield an enormous amount of important information, including the chromatic polynomial, the flow polynomial, the all-terminal reliability polynomial,
the weight enumerator of a linear code, and numbers of spanning trees, forests, spanning subgraphs, and acyclic orientations.  (See \cite{ellis-monaghan-moffatt2022}.)

Dominic used to promote research on the Tutte polynomial using diagrams of the ``Tutte plane'', drawn at coffee tables or displayed on slides in talks.  (Graham Farr recalls one from 1983, but there were surely others before then.)  These diagrams show the real $xy$-plane together with various points and curves where $T(G;x,y)$ has interesting interpretations.  We give such a diagram in Figure~\ref{fig:tutte-plane}.  Diagrams of this type may also be found in \citelist{\cite{welsh-93b}*{Fig.~8.2} \cite{farr-07}*{Fig.~3.1}}.  The diagram illustrates, for example, that
the chromatic polynomial $P(G;\lambda)$ of a graph $G=(V,E)$ can be found from $T(G;x,y)$ along the line $y=0$ by appropriate substitutions and simple algebra.
In this case, the actual relationship --- due to Whitney \cite{whitney1932} and using (\ref{eq:tutte-from-whitney}) --- is
\[
P(G;\lambda) = (-1)^{r(E)}\lambda^{|V|-r(E)}T(G; 1-\lambda, 0).
\]
The hyperbolae $H_q=\{(x,y)\mid (x-1)(y-1)=q\}$ are significant: appropriate parameterisations give the partition function of the $q$-state Potts model, and, for $q=2$, the closely related Ising model and (at the more general level of binary matroids) the weight enumerator of a linear code over GF(2).  On the hyperbola $H_1$, the Tutte polynomial simplifies to $(x-1)^{r(E)}y^{|E|}$, which is easy to compute but degenerate and uninformative; the only matroid properties it contains are the rank and size of the ground set.  Not all pertinent hyperbolae are of the form $H_q$;
the hyperbola $\{(x,y)\mid xy=1\}$ carries the Jones polynomial of an alternating link (and the ``degenerate hyperbola'' $xy=0$ carries the chromatic and flow polynomials).  Our diagram focuses on graph-theoretic interpretations, with some matroid interpretations mentioned in parentheses in the legend.  Further details are in \cite{ellis-monaghan-moffatt2022,welsh-93b}.

\begin{figure}
\begin{center}
\begin{tikzpicture}
\begin{axis}[
legend style={font=\small},
legend cell align=left,
legend style={at={(0.5,-0.05)},anchor=north},
legend columns=7, 
transpose legend,
unit vector ratio = {1 1},
width=14cm,
height=14cm,
axis x line=center,
  axis y line=center,
  xtick={-4,-3,...,4},
  x tick label style={xshift={1.3mm}},
  xticklabels={,,,,,1,2,,},
  ytick={-4,-3,...,4},
  yticklabels={,,,,,1,2,,},
  xlabel={$x$},
  ylabel={$y$},
  xlabel style={below right},
  ylabel style={above left},
  xmin=-5.8,
  xmax=5.8,
  ymin=-5.8,
  ymax=5.8,
    samples=50,
    axis line style={latex-latex},
]
  \addplot[dashed,domain=-5.5:5.5,forget plot] (x,1) ;
  \addplot[dashed,domain=-5.5:5.5,forget plot] (1,x) ;

  \addplot[greenish,ultra thick,domain=-5.5:5.5] (x,0) ;
  \addlegendentry{$y=0$: chromatic};

  \addplot[olive,ultra thick,domain=-5.5:5.5] (0,x) ;
  \addlegendentry{$x=0$: flow};

  \addplot[purple,ultra thick,domain=-5.5:5.5] (1,x) ;
  \addlegendentry{$x=1$: reliability};

 \addplot[black, ultra thick, domain=1.222:5.5] {1+1/(x-1)} ;   
  \addlegendentry{$H_1$: easy curve};
  \addplot[black, ultra thick, domain=-5.5:0.846,forget plot] {1+1/(x-1)} ;   
  
  \addplot[blue, ultra thick, domain=1.444:5.5] {1+2/(x-1)} ;   
  \addlegendentry{$H_2$: Ising (also, weight enum.)~~~~~~};
  \addplot[blue, ultra thick, domain=-5.5:0.692,forget plot] {1+2/(x-1)} ;   
  
  \addplot[violet, ultra thick, domain=2.111:5.5] {1+5/(x-1)} ;   
  \addlegendentry{$H_q$: $q$-state Potts~~~};
  \addplot[violet, ultra thick, domain=-5.5:0.231,forget plot] {1+5/(x-1)} ;   
  
  \addplot[ochre, ultra thick, domain=0.1818:5.5] {1/x} ;   
  \addlegendentry{$xy=1$: Jones};
  \addplot[ochre, ultra thick, domain=-5.5:-0.1818,forget plot] {1/x} ;   

  \addplot[black,ultra thick,fill=purple,only marks,mark size=1.4mm] coordinates {(1,1)};
  \addlegendentry{\# spanning trees (bases)};
  
  \addplot[red,only marks,mark size=1.5mm,mark=triangle*] coordinates {(2,1)};
  \addlegendentry{\# forests (indep.\ sets)};

  \addplot[red,fill,only marks,mark size=1.5mm,mark=triangle*,every mark/.append style={rotate=180}] coordinates {(1,2)};
  \addlegendentry{\# spanning sets};

  \addplot[greenish,fill,only marks,mark size=1.5mm,mark=pentagon*] coordinates {(2,0)};
  \addlegendentry{\# acyclic orientations};

  \addplot[olive,fill,only marks,mark size=1.5mm,mark=pentagon*,every mark/.append style={rotate=180}] coordinates {(0,2)};
  \addlegendentry{\# totally cyclic orientations};

  \addplot[black,ultra thick,fill=greenish,only marks,mark size=1.4mm,mark=square*] coordinates {(-1,0)};
  \addlegendentry{\# 2-colourings};

  \addplot[black,ultra thick,fill=olive,only marks,mark size=1.4mm,mark=square*,every mark/.append style={rotate=45}] coordinates {(0,-1)};
  \addlegendentry{\# 2-colourings of dual};

 \addplot[black,ultra thick,fill=ochre,only marks,mark=oplus*,mark size=1.4mm,forget plot] coordinates {(-1,-1)} ;

\end{axis}
\end{tikzpicture}
\end{center}
\caption{The Tutte plane.}
\label{fig:tutte-plane}
\end{figure}

Much of the information in the Tutte polynomial is \#P-hard
to compute; \#P-hardness is a stronger form of NP-hardness for counting problems.  For example, counting 3-colourings of a graph is \#P-hard \cite{linial1986}.  So computing the entire Tutte polynomial is \#P-hard too.

In the late 1980s, Dominic combined his interests in matroid theory and complexity to investigate the complexity of evaluating the Tutte polynomial at specific points.  Work with Dirk Vertigan (DPhil, 1991) and Fran\c{c}ois Jaeger, which comprised part of the latter's DPhil thesis, resulted in a complete classification of points $(x,y)$ according to whether computation of $T(G;x,y)$ is polynomial time or \#P-hard
\cite{JVW-90}.  This paper was published in 1990 and yielded new complexity results for some other combinatorial polynomials that can be obtained from the Tutte polynomial by algebraic substitutions, including the Jones polynomial from knot theory and the partition functions of the Ising and Potts models from statistical physics.  For $x,y\in\mathbb{R}$, Figure~\ref{fig:tutte-plane} illustrates the result.  The points $(x,y)$
where $T(G;x,y)$ is polynomial-time computable are the four points indicated
by symbols outlined in black and all points on the black hyperbola $H_1$;
evaluation at all other points is \#P-hard.  The result of Jaeger, Vertigan and Welsh \cite{JVW-90} covers all
complex $x,y$, giving four further polynomial-time computable points
with all other cases being \#P-hard.  Soon afterwards, Dirk Vertigan and Dominic gave a similar complexity classification for the restriction to bipartite planar graphs \cite{vertigan-welsh1992}.
For his doctoral research on the complexity of Tutte polynomial evaluation at specific points, Vertigan was awarded the Senior Mathematical Prize at Oxford in 1990 \cite{senior-mathematical-prize-1990}.
This was an annual prize for ``the dissertation of greatest merit on any subject of Pure or Applied Mathematics'' \cite{senior-mathematical-prize-1990-regs}.

This work opened up a rich new vein of research in which complexity results were proved for various properties (including evaluations and coefficients) of various polynomials for various classes of graphs and matroids.  Some of these results pertained to exact computation of the polynomials, others to approximating or bounding them.  The main tool for approximation was a type of efficient randomised algorithm --- a fully polynomial randomised approximation scheme (FPRAS) --- and the main approach was to develop a Markov Chain Monte Carlo (MCMC) algorithm, which uses a Markov chain among the structures of interest to sample from them and hence estimate their numbers.  A key technical challenge is to ensure that the Markov chain converges quickly enough (is ``rapidly mixing''), and issues of this type led Dominic to some significant combinatorial problems and conjectures.

He wrote a book \textit{Complexity: Knots, Colourings and Counting} \cite{welsh-93b} which for many years was the leading reference on Tutte polynomials, their complexity and their links with other areas.  The book is still widely used.  His many articles on the topic included updated surveys which were very accessible and influential (see \cite{farr-07}).  His significant survey articles include \cite{welsh1999}, on the Tutte polynomial, and \cite{welsh1997b} on approximate counting in general.

This stream of research also includes the MSc of Laura Ch\'avez Lomelí (1994) and the DPhils of James Annan (1994), Steve Noble (1997), Eric Bartels (1997), and Magnus Bordewich (2003).  Annan proved that computing any specific coefficient of the Tutte polynomial is \#P-complete \cite{annan-95} and gave, for dense graphs, an FPRAS for the value of the Tutte polynomial at $(x,1)$ where $x \ge 1$, which includes counting forests ($x=2$) \cite{annan-94}.  Dominic later collaborated with Noga Alon and Alan Frieze to extend this result to evaluations at any $(x,y)$ with $x,y \ge 1$ \cite{AFW-95}.  Whether a randomised approximation of this type can be found for all graphs remains an open question, though Bordewich was able to find such approximations for certain types of sparse graphs \cite{bordewich-04}.  Noble proved that the Tutte polynomial can be evaluated in polynomial time for graphs of bounded tree-width \cite{noble-98}, and, in \cite{noble-07}, he provided a Jaeger--Vertigan--Welsh-style complexity classification for evaluations of an analogue of the Tutte polynomial for 2-polymatroids introduced by Oxley and Whittle in \cites{oxley-whittle1993a,oxley-whittle1993b}.

Dominic collaborated with Bordewich, Freedman, and Lov\'asz on an important paper \cite{BFLW-05} showing that an additive approximation (which is weaker than an FPRAS) to a certain Tutte polynomial evaluation (related to the Jones polynomial) is sufficient to capture the power of quantum computation, extending earlier work of Freedman, Kitaev, Larson, and Wang \cite{FKLW-03}.  See Ken Regan's tribute \cite{regan2024} for further discussion of this.

Dominic used MCMC algorithms in some other novel ways.  He developed a Markov chain method for choosing a planar graph uniformly at random amongst those with a fixed vertex set of size $n$ (with Denise and Vasconcellos) \cite{DVW-96}.  He later used this idea with McDiarmid and Steger (2005) and discovered the surprising fact that the probability that a random planar graph is connected tends to neither 0 nor 1 as $n$ tends to infinity \cite{MSW-05}.  Their theory is developed further in \cite{MSW-2006}.

Dominic also explored some purely combinatorial properties of Tutte--Whitney polynomials and found interesting generalisations, evaluations and relationships.  For example, he proved that the expected values of Tutte polynomial evaluations for a random subgraph of a graph are themselves evaluations of the Tutte polynomial of the given graph \cite{welsh-96}.  With Steve Noble, he introduced a more general family of graph polynomials inspired by some knot invariants \cite{noble-welsh1999}.  With Geoff Whittle, he introduced a different generalisation inspired by hyperplane arrangements, which turned out to count certain channel assignments in communication theory \cite{welsh-whittle1999}.  He took a particular interest in the way Tutte--Whitney polynomials linked different fields of mathematics: combinatorics, statistical physics, network reliability, coding theory, and knot theory, and wrote often on these links.  He fostered connections with researchers in these diverse fields.

In a 1995 paper, Dominic and his then-student Eric Bartels defined the mean colour number of an $n$-vertex graph to be the expected number of colours actually used in a random proper $n$-colouring of the graph, and conjectured that it achieves its minimum when the graph is empty \cite{bartels-welsh-95}.  Dominic called this the Shameful Conjecture because so many researchers had tried but failed to prove it in spite of it seeming so obvious.  It finally succumbed to Fengming Dong \cite{dong-00}.  For an account of its history, see \cite{royle-15}.  The Bartels--Welsh paper also proposed a stronger conjecture, that the mean colour number never decreases when an edge is added to a graph, but this was disproved by Michele Mosca \cite{mosca-98} who later completed a DPhil (1999) in quantum computing under the joint supervision of Dominic and his colleague Artur Ekert.  These conjectures have the character of correlation inequalities and share their combination of plausibility and difficulty.  They were not the last correlation inequalities to be conjectured by Dominic and acquire fame (or notoriety!).

His work on other aspects of the Tutte polynomial included supervision of the DPhils of Criel Merino (2000), Koko Kayibi (2002), and Andrew Goodall (2004).  Merino linked certain evaluations of the Tutte polynomial to numbers of critical configurations in a chip-firing game on graphs \cite{merino-97}.  He also did a computational study of the asymptotic behaviour of the number of acyclic orientations in square-lattice graphs (another link to some of Dominic’s earliest interests) \cite{merino-welsh-99}.  This led him and Dominic to conjecture that the number of spanning trees of a graph is bounded above by whichever is larger of the number of acyclic orientations and the number of totally cyclic orientations \cite{merino-welsh-99}*{Conj.~7.1}.  This too may be viewed as a correlation inequality.  This conjecture is known as the Merino--Welsh Conjecture and remains open.  It has proved to be very attractive and very difficult.  Gordon Royle has posted a good introduction at \cite{royle-12}.  For general matroids, counterexamples have recently been given by Beke, Cs\'aji, Csikv\'ari, and Pituk \cite{BCCP-24}.

Dominic co-organised the first workshop on the Tutte polynomial, at Barcelona in 2001.  It became the basis of a special issue of \textit{Advances in Applied Mathematics} in 2004, which he co-edited \cite{kung-noy-welsh2004}.  The workshop was a milestone in the building of a strong community of researchers on this topic.

He retained his early interest in the complexity of matroid computations in general, and his last DPhil student, Dillon Mayhew (2005), took up this theme, this time considering matroids to be represented (for input purposes) not by an oracle (as Robinson and Welsh had done \cite{robinson-welsh-80}) but by a complete list of all the sets required for one of the axiomatisations, that is, all independent sets, or all bases, or all circuits, and so on \cite{mayhew-08}.

Research on other matroid topics continued through this fourth phase, including through his supervision of Irasema Sarmiento (DPhil, 1998) and (with Colin McDiarmid) Rhiannon Hall (2005).  Once Dominic retired in 2005, he seemed to spend more time on matroid theory, for example in studying matroid enumeration and properties of random matroids with several collaborators in papers published in the period 2011--2013 \cites{MNWW-11, LOSW-13, mayhew-welsh-13}.

His very last paper, on counting phylogenetic networks (with Colin McDiarmid and Charles Semple), showed that he retained his sense of mathematical adventure in tackling new problems \cite{mcdiarmid-semple-welsh2015}.

\section{The complete mathematician}
\label{sec:complete-mathematician}

Dominic made lasting contributions to mathematics.  His research style was characterised by a deep sense of mathematical beauty and astute judgement as to what problems or topics were worth pursuing.  In framing theorems and writing about them, he was good at drawing out and highlighting what was really important.  Dominic had a natural gift for linking disparate fields, a real sense of intellectual adventure, and a youthful readiness to move into new territory.  He seemed impelled by a kind of level-headed urgency to make the most of opportunities and push mathematics forward.

Dominic’s influence goes well beyond what can be accounted for by his formal contributions to building the intellectual structure of mathematics --- his definitions, theorems and proofs.  He also advanced mathematics by his development of its people, community, and culture.  Not only was he eminently successful at mathematics, he was a complete mathematician.

This is evident firstly in his writing, through his many survey papers and textbooks.  Moreover, many of his research articles contain valuable expositions of the state of particular research topics and their relationships with other fields.  These played a significant and easily underestimated role in the development of the fields he worked in.  His writings also show an interest in history, an awareness of emerging unpublished work, and liberal acknowledgements of the contributions of others.

A striking characteristic of Dominic’s expository power --- in writing, presenting, and conversation --- was his extraordinary ability to move rapidly from introducing an area to stating enticing research problems in that area.  Mathematical depth was made accessible and shared.  His diagrams of the Tutte plane and of the world of complexity theory were ubiquitous in his talks and they were memorable for the compact way they encapsulated so much information.  Our Figures \ref{fig:complexity} and \ref{fig:tutte-plane} are inspired by those diagrams.

His advancement of the human side of mathematics was aided greatly by his character and social attributes.  He had a remarkable ability to get on well with people, indeed his personality had a magnetic quality.  He was equal to any social situation and always pleasant, considerate, and courteous, often smiling or laughing.  He was very hospitable and, together with his wife Bridget, he would welcome students and colleagues to the home he shared with her and their sons.  He was a fine athlete and was fiercely competitive, challenging all comers variously in real tennis, table tennis, squash, croquet, boules and speed chess.  All such contests had laughter associated with them, though sometimes this occurred long after the event.

Dominic was at once outgoing and sensitive, confident and careful, firm and gentle, purposeful and relaxed, serious and light-hearted, dignified and informal.  These somewhat-contrasting qualities would often play out together as he talked.  He was quick to find humour in situations and freely shared anecdotes and observations.  Anyone interacting with him would be confident of his attention and goodwill.  He built strong and enduring connections.  These days, one might say that he was a central node in the social network of his field.  He used that position generously in sharing news, ideas, and problems.  His informal interactions with other researchers were very influential and have helped shape the areas he worked in.

Dominic led a strong combinatorics group at the Mathematical Institute in Oxford, with colleagues Peter Cameron (who recently reflected on his influence at \cite{cameron-23}), Aubrey Ingleton, Colin McDiarmid, and their students.  He was a popular and effective teacher of undergraduates.  Under the Oxford tutorial system, he worked most closely with undergraduates at Merton College where he was a major contributor to his College’s very strong results in mathematics over many years.  Undergraduates taught by Dominic at Merton included, chronologically, Andrew Wiles, Dugald Macpherson, Peter Kronheimer, and Dominic Joyce, who between them have gone on to win more than twenty major mathematical awards including the Abel Prize to Wiles.

He was very efficient and well organised, with a deep sense of duty.  So it is not surprising that he took on leadership roles.  He edited several journals, organised some pivotal conferences, led the British Combinatorial Committee, and held senior positions at Oxford, to the great benefit of all his academic communities.  In 1985, as second chair of the British Combinatorial Committee, Dominic gave a memorable after-dinner address at a reception hosted by the Lord Provost of Glasgow at the Glasgow City Chambers.  Dominic's theme was the attraction of doing mathematics and, echoing G.H. Hardy, he said that, for him, the joy came from discovering beautiful patterns.  All who are acquainted with Dominic's work, even superficially, will recognise how this search for such discoveries permeates his research.

Dominic was a very effective supervisor of research students.  This drew not only on his mathematical knowledge, taste, judgement, and connections, but also on his personality, character, and people skills.  He was flexible in his approach and adept at finding a productive mix of patience, firmness, encouragement, plain speaking, and inspiration.  He was generous in the freedom he gave his students while being as supportive as needed.  Time and again he brought out the best in his diverse range of students, inspiring enduring appreciation and affection.  An academic family tree for Dominic, listing his DPhil graduates (28 as main supervisor) and many of his academic descendants, is maintained by David Wood (Graham Farr’s first PhD graduate, 2000) at \cite{wood-1}.  This is much more comprehensive than the one at the Mathematics Genealogy Project \cite{math-genealogy-project}.  David also maintains the academic family tree of John Hammersley \cite{wood-2}, of which Dominic’s tree is the dominant subtree.

Dominic Welsh's enduring legacy includes a significant and eclectic body of research, a rich library of influential expository writing, and a large and notable community of students, collaborators, and colleagues.  He will be remembered with admiration, affection, and deep gratitude.  He has greatly enriched discrete mathematics as both a research field and a human endeavour.   \\

\noindent\textbf{Acknowledgements}

We are grateful for comments and suggestions from Peter Cameron, Keith Edwards, Geoffrey Grimmett, Colin McDiarmid, Steven Noble, Ken Regan, Charles Semple, Geoff Whittle, and David Wood.  Thanks also to Tony Huynh and the Matroid Union for their encouragement and assistance with the first version of this article \cite{farr-mayhew-oxley2024}.  We are grateful to Dominic's widow, Bridget Welsh, for kindly giving permission to the Bodleian Library to scan Dominic's DPhil thesis \cite{welsh-64} and make it available online via the Library's website.  We thank the Monash University Library for funding that process for the thesis and for the assistance of their Resource Sharing Team.  We also thank Ruth Anderson (Deputy Editor, University of Oxford Gazette) and Catherine McIlwaine (Archivist, Special Collections, Bodleian Libraries) for information on \cite{senior-mathematical-prize-1990,senior-mathematical-prize-1990-regs}.


\begin{bibdiv}

\begin{biblist}

\bib{fields}{webpage}{
   title={Fields Medals 2022},
   url={https://www.mathunion.org/imu-awards/fields-medal/fields-medals-2022},
   accessdate={2024-04-20},   
}

\bib{math-genealogy-project}{webpage}{
    title={Mathematics Genealogy Project},
    url={https://www.mathgenealogy.org},
    note={North Dakota State University Department of Mathematics},
    accessdate={2024-04-20},
}

 
\bib{senior-mathematical-prize-1990}{article}{
   title={Senior Mathematical Prize 1990},
   journal={Oxford University Gazette},
   volume={CXX},
   date={28 June 1990},   
   number={4185},
   pages={1001},
}

\bib{senior-mathematical-prize-1990-regs}{misc}{
   title={Statutes, Decrees and Regulations of the University of Oxford},
   note={Ch. IX, section I.174.1--2},
   date={1989},
}

\bib{AHK-18}{article}{
   author={Adiprasito, Karim},
   author={Huh, June},
   author={Katz, Eric},
   title={Hodge theory for combinatorial geometries},
   journal={Ann. of Math. (2)},
   volume={188},
   date={2018},
   number={2},
   pages={381--452},
}

\bib{AFW-95}{article}{
   author={Alon, Noga},
   author={Frieze, Alan},
   author={Welsh, Dominic},
   title={Polynomial time randomized approximation schemes for
   Tutte-Gr\"{o}thendieck invariants: the dense case},
   journal={Random Structures Algorithms},
   volume={6},
   date={1995},
   number={4},
   pages={459--478},
}

\bib{annan-94}{article}{
   author={Annan, J. D.},
   title={A randomised approximation algorithm for counting the number of
   forests in dense graphs},
   journal={Combin. Probab. Comput.},
   volume={3},
   date={1994},
   number={3},
   pages={273--283},
}

\bib{annan-95}{article}{
   author={Annan, J. D.},
   title={The complexities of the coefficients of the Tutte polynomial},
   note={Combinatorial optimization 1992 (CO92) (Oxford)},
   journal={Discrete Appl. Math.},
   volume={57},
   date={1995},
   number={2-3},
   pages={93--103},
}

\bib{BPP-14}{article}{
   author={Bansal, N.},
   author={Pendavingh, R. A.},
   author={van der Pol, J. G.},
   title={An entropy argument for counting matroids},
   journal={J. Combin. Theory Ser. B},
   volume={109},
   date={2014},
   pages={258--262},
}

\bib{BPP-15}{article}{
   author={Bansal, Nikhil},
   author={Pendavingh, Rudi A.},
   author={van der Pol, Jorn G.},
   title={On the number of matroids},
   journal={Combinatorica},
   volume={35},
   date={2015},
   number={3},
   pages={253--277},
}

\bib{bartels-welsh-95}{article}{
   author={Bartels, J. Eric},
   author={Welsh, Dominic J. A.},
   title={The Markov chain of colourings},
   conference={
      title={Integer programming and combinatorial optimization},
      address={Copenhagen},
      date={1995},
   },
   book={
      series={Lecture Notes in Comput. Sci.},
      volume={920},
      publisher={Springer, Berlin},
   },
   date={1995},
   pages={373--387},
}

\bib{BCCP-24}{article}{
   author={Beke, Csongor},
   author={Cs\'{a}ji, Gergely K\'{a}l},
   author={Csikv\'{a}ri, P\'{e}ter},
   author={Pituk, S\'{a}ra},
   title={The Merino--Welsh conjecture is false for matroids},
   journal={Adv. Math.},
   volume={446},
   date={2024},
   pages={Paper No. 109674},
}

\bib{biggs1974}{book}{
   author={Biggs, Norman},
   title={Algebraic graph theory},
   series={Cambridge Tracts in Mathematics},
   volume={No. 67},
   publisher={Cambridge University Press, London},
   date={1974},
   pages={vii+170},
}

\bib{bondy-welsh-66}{article}{
   author={Bondy, J. A.},
   author={Welsh, D. J. A.},
   title={A note on the monomer dimer problem},
   journal={Proc. Cambridge Philos. Soc.},
   volume={62},
   date={1966},
   pages={503--505},
}

\bib{bondy-welsh-71}{article}{
   author={Bondy, J. A.},
   author={Welsh, D. J. A.},
   title={Some results on transversal matroids and constructions for
   identically self-dual matroids},
   journal={Quart. J. Math. Oxford Ser. (2)},
   volume={22},
   date={1971},
   pages={435--451},
}

\bib{bordewich-04}{article}{
   author={Bordewich, Magnus},
   title={Approximating the number of acyclic orientations for a class of
   sparse graphs},
   journal={Combin. Probab. Comput.},
   volume={13},
   date={2004},
   number={1},
   pages={1--16},
}

\bib{BFLW-05}{article}{
   author={Bordewich, M.},
   author={Freedman, M.},
   author={Lov\'{a}sz, L.},
   author={Welsh, D.},
   title={Approximate counting and quantum computation},
   journal={Combin. Probab. Comput.},
   volume={14},
   date={2005},
   number={5-6},
   pages={737--754},
}

\bib{broadbent-hammersley1957}{article}{
   author={Broadbent, S. R.},
   author={Hammersley, J. M.},
   title={Percolation processes. I. Crystals and mazes},
   journal={Proc. Cambridge Philos. Soc.},
   volume={53},
   date={1957},
   number={5-6},
   pages={629--641},
}

\bib{brylawski-82}{article}{
   author={Brylawski, Thomas},
   title={The Tutte polynomial. I. General theory},
   conference={
      title={Matroid theory and its applications},
   },
   book={
      publisher={Liguori, Naples},
   },
   date={1982},
   pages={125--275},
}

\bib{cameron-23}{webpage}{
    author={Cameron, P. J.},
    title={Dominic Welsh},
    date={2023-12-04},
    url={https://cameroncounts.wordpress.com/2023/12/04/dominic-welsh/},
    accessdate={2024-04-20},
}

\bib{chavezlomeli-welsh1996}{article}{
   author={Ch\'avez Lomel\'\i, Laura},
   author={Welsh, Dominic},
   title={Randomised approximation of the number of bases},
   conference={
      title={Proc. AMS-IMS-SIAM Joint Summer Research Conf.},
      address={University of Washington, Seattle, WA},
      date={July 2–6, 1995},
   },
   book={
      editor={Bonin, Joseph E.},
      editor={Oxley, James G.},
      editor={Servatius, Brigitte},
      title={Matroid theory},
      series={Contemp. Math.},
      volume={197},
      publisher={Amer. Math. Soc., Providence, RI},
   },
   date={1996},
   pages={371--376},
}

\bib{crapo-schmitt-05}{article}{
   author={Crapo, Henry},
   author={Schmitt, William},
   title={A unique factorization theorem for matroids},
   journal={J. Combin. Theory Ser. A},
   volume={112},
   date={2005},
   number={2},
   pages={222--249},
}

\bib{DVW-96}{article}{
   author={Denise, Alain},
   author={Vasconcellos, Marcio},
   author={Welsh, Dominic J. A.},
   title={The random planar graph},
   note={Festschrift for C. St. J. A. Nash-Williams},
   journal={Congr. Numer.},
   volume={113},
   date={1996},
   pages={61--79},
}

\bib{desousa-welsh-72}{article}{
   author={de Sousa, J.},
   author={Welsh, D. J. A.},
   title={A characterisation of binary transversal structures},
   journal={J. Math. Anal. Appl.},
   volume={40},
   date={1972},
   pages={55--59},
}

\bib{dong-00}{article}{
   author={Dong, F. M.},
   title={Proof of a chromatic polynomial conjecture},
   journal={J. Combin. Theory Ser. B},
   volume={78},
   date={2000},
   number={1},
   pages={35--44},
}

\bib{donnelly-welsh-83}{article}{
   author={Donnelly, Peter},
   author={Welsh, Dominic},
   title={The antivoter problem: random $2$-colourings of graphs},
   conference={
      title={Graph theory and combinatorics},
      address={Cambridge},
      date={1983},
   },
   book={
      publisher={Academic Press, London},
   },
   date={1984},
   pages={133--144},
}

\bib{edwards-85}{article}{
   author={Edwards, Keith},
   title={Counting self-avoiding walks in a bounded region},
   journal={Ars Combin.},
   volume={20},
   date={1985},
   pages={271--281},
}

\bib{edwards1986}{article}{
   author={Edwards, Keith},
   title={The complexity of colouring problems on dense graphs},
   journal={Theoret. Comput. Sci.},
   volume={43},
   date={1986},
   number={2--3},
   pages={337--343},
}

\bib{edwards1992}{article}{
   author={Edwards, Keith},
   title={The complexity of some graph colouring problems},
   journal={Discrete Appl. Math.},
   volume={36},
   date={1992},
   number={2},
   pages={131--140},
}

\bib{ellis-monaghan-moffatt2022}{collection}{
   title={Handbook on the Tutte Polynomial and Related Topics},
   editor={Ellis-Monaghan, Joanna},
   editor={Moffatt, Iain},
   publisher={Chapman and Hall/CRC Press},
   date={2022},
}

\bib{farr-07}{article}{
   author={Farr, Graham E.},
   title={Tutte-Whitney polynomials: some history and generalizations},
   book={
      title={Combinatorics, complexity, and chance},
      editor={Grimmett, Geoffrey},
      editor={McDiarmid, Colin},
      series={Oxford Lecture Ser. Math. Appl.},
      volume={34},
      publisher={Oxford Univ. Press, Oxford},
   },
   date={2007},
   pages={28--52},
}

\bib{farr-mayhew-oxley2024}{webpage}{
    author={Farr, G.},
    author={Mayhew, D.},
    author={Oxley, J.},
    title={Dominic Welsh: his work and influence},
    date={2024-02-12},
    url={http://matroidunion.org/?p=5304},
    accessdate={2024-04-20},
}

\bib{FKLW-03}{article}{
   author={Freedman, Michael H.},
   author={Kitaev, Alexei},
   author={Larsen, Michael J.},
   author={Wang, Zhenghan},
   title={Topological quantum computation},
   note={Mathematical challenges of the 21st century (Los Angeles, CA,
   2000)},
   journal={Bull. Amer. Math. Soc. (N.S.)},
   volume={40},
   date={2003},
   number={1},
   pages={31--38},
}

\bib{garey-johnson1979}{book}{
   author={Garey, Michael R.},
   author={Johnson, David S.},
   title={Computers and Intractability: A Guide to the Theory of NP-Completeness},
   publisher={W. H. Freeman, San Francisco},
   date={1979},
   pages={x+338},
}

\bib{grimmett-89}{book}{
   author={Grimmett, Geoffrey},
   title={Percolation},
   publisher={Springer-Verlag, New York},
   date={1989},
   pages={xii+296},
}

\bib{grimmett2024a}{article}{
   author={Grimmett, Geoffrey R.},
   title={Dominic Welsh (1938–2023)},
   journal={Bull. Lond. Math. Soc.},
   status={to appear},
   note={Preprint:  \url{https://arxiv.org/pdf/2404.13942}},
}

\bib{grimmett2024b}{webpage}{
   author={Grimmett, Geoffrey R.},
   title={Publications of Dominic Welsh},
   url={https://www.statslab.cam.ac.uk/~grg/papers/welsh-bib.pdf},
   note={Accompanies \cite{grimmett2024a}},
}


\bib{grimmett-mcdiarmid2007}{collection}{
   title={Combinatorics, complexity, and chance},
   series={Oxford Lecture Series in Mathematics and its Applications},
   volume={34},
   editor={Grimmett, Geoffrey},
   editor={McDiarmid, Colin},
   note={A tribute to Dominic Welsh},
   publisher={Oxford University Press, Oxford},
   date={2007},
   pages={x+310},
}

\bib{grimmett-welsh-86}{book}{
   author={Grimmett, Geoffrey},
   author={Welsh, Dominic},
   title={Probability: an introduction},
   series={Oxford Science Publications},
   publisher={The Clarendon Press, Oxford University Press, New York},
   date={1986},
}

\bib{grimmett-welsh-14}{book}{
   author={Grimmett, Geoffrey},
   author={Welsh, Dominic},
   title={Probability---an introduction},
   edition={2},
   publisher={Oxford University Press, Oxford},
   date={2014},
}

\bib{hammersley-57a}{article}{
   author={Hammersley, J. M.},
   title={Percolation processes. II. The connective constant},
   journal={Proc. Cambridge Philos. Soc.},
   volume={53},
   date={1957},
   pages={642--645},
}

\bib{hammersley-57b}{article}{
   author={Hammersley, J. M.},
   title={Percolation processes: Lower bounds for the critical probability},
   journal={Ann. Math. Statist.},
   volume={28},
   date={1957},
   pages={790--795},
}

\bib{hammersley-morton1954}{article}{
   author={Hammersley, J. M.},
   author={Morton, K. W.},
   title={Poor man's Monte Carlo},
   journal={J. Roy. Statist. Soc. Ser. B},
   volume={16},
   date={1954},
   pages={23--38; discussion 61--75},
}

\bib{hammersley-welsh1962}{article}{
   author={Hammersley, J. M.},
   author={Welsh, D. J. A.},
   title={Further results on the rate of convergence to the connective constant of the hypercubical lattice},
   journal={Quart. J. Math. Oxford Ser. (2)},
   volume={13},
   date={1962},
   pages={108--110},
}

\bib{hammersley-welsh1965}{article}{
   author={Hammersley, J. M.},
   author={Welsh, D. J. A.},
   title={First-passage percolation, subadditive processes, stochastic networks, and generalized renewal theory},
   conference={
      title={Proc. Internat. Res. Semin., Statist. Lab., Univ. California, Berkeley, Calif., 1963},
   },
   book={
      editor={Jerzy Neyman},
      editor={Lucien M. LeCam},
      title={Bernoulli 1713, Bayes 1763, Laplace 1813: Anniversary volume},
      publisher={Springer-Verlag, New York},
   },
   date={1965},
   pages={61--110},
}

\bib{hammersley-welsh1980}{article}{
   author={Hammersley, J. M.},
   author={Welsh, D. J. A.},
   title={Percolation theory and its ramifications},
   journal={Contemp. Phys.},
   volume={21},
   number={6},
   date={1980},
   pages={593-605},
}


\bib{harary-welsh-68}{article}{
   author={Harary, Frank},
   author={Welsh, Dominic},
   title={Matroids versus graphs},
   conference={
      title={The Many Facets of Graph Theory},
      address={Proc. Conf., Western Mich. Univ., Kalamazoo, Mich.},
      date={1968},
   },
   book={
      series={Lecture Notes in Math.},
      volume={No. 110},
      publisher={Springer, Berlin-New York},
   },
   date={1969},
   pages={155--170},
}

\bib{heron-72}{article}{
   author={Heron, A. P.},
   title={Matroid polynomials},
   conference={
      title={Combinatorics},
      address={Proc. Conf. Combinatorial Math., Math. Inst., Oxford},
      date={1972},
   },
   book={
      publisher={Inst. Math. Appl., Southend-on-Sea},
   },
   date={1972},
   pages={164--202},
}



\bib{JVW-90}{article}{
   author={Jaeger, F.},
   author={Vertigan, D. L.},
   author={Welsh, D. J. A.},
   title={On the computational complexity of the Jones and Tutte
   polynomials},
   journal={Math. Proc. Cambridge Philos. Soc.},
   volume={108},
   date={1990},
   number={1},
   pages={35--53},
}

\bib{kesten-80}{article}{
   author={Kesten, Harry},
   title={The critical probability of bond percolation on the square lattice
   equals ${1\over 2}$},
   journal={Comm. Math. Phys.},
   volume={74},
   date={1980},
   number={1},
   pages={41--59},
}


\bib{kung-noy-welsh2004}{misc}{
   editor={Kung, J. P. S.},
   editor={Noy, M.},
   editor={Welsh, D. J. A.},
   title={Special Issue on the Tutte Polynomial, {\rm Adv. in Appl. Math. \textbf{32} (Jan.--Feb. 2004), nos. 1--2}},
}

\bib{lemos-04}{article}{
   author={Lemos, Manoel},
   title={On the number of non-isomorphic matroids},
   journal={Adv. in Appl. Math.},
   volume={33},
   date={2004},
   number={4},
   pages={733--746},
}

\bib{lenz-13}{arXiv}{
  author={Lenz, Matthias},
  title={Matroids and log-concavity},
  date={2013},
  eprint={1106.2944},
  archiveprefix={arXiv},
  primaryclass={math.CO},
}

\bib{linial1986}{article}{
   author={Linial, Nathan},
   title={Hard enumeration problems in geometry and combinatorics},
   journal={SIAM J. Alg. Disc. Meth.},
   volume={7},
   date={1986},
   number={2},
   pages={331--335},
}

\bib{lipton-regan}{webpage}{
    author={Lipton, R.},
    author={Regan, K.},
    title={G\"{o}del’s Lost Letter and P=NP},
    url={https://rjlipton.wpcomstaging.com},
    accessdate={2024-04-20},
}

\bib{LOSW-13}{article}{
   author={Lowrance, Lisa},
   author={Oxley, James},
   author={Semple, Charles},
   author={Welsh, Dominic},
   title={On properties of almost all matroids},
   journal={Adv. in Appl. Math.},
   volume={50},
   date={2013},
   number={1},
   pages={115--124},
}

\bib{mcdiarmid-semple-welsh2015}{article}{
   author={McDiarmid, Colin},
   author={Semple, Charles},
   author={Welsh, Dominic J. A.},
   title={Counting phylogenetic networks},
   journal={Ann. Comb.},
   volume={19},
   date={2015},
   number={1},
   pages={205--224},
}

\bib{MSW-05}{article}{
   author={McDiarmid, Colin},
   author={Steger, Angelika},
   author={Welsh, Dominic J. A.},
   title={Random planar graphs},
   journal={J. Combin. Theory Ser. B},
   volume={93},
   date={2005},
   number={2},
   pages={187--205},
}


\bib{MSW-2006}{article}{
   author={McDiarmid, Colin},
   author={Steger, Angelika},
   author={Welsh, Dominic J. A.},
   title={Random graphs from planar and other addable classes},
   book={
      title={Topics in discrete mathematics},
      editor={Klazar, M.},
      editor={Kratochv\'\i l, Jan},
      editor={Loebl, Martin},
      editor={Matou\v sek, Ji\v r\'\i},
      editor={Thomas, Robin},
      editor={Valtr, Pavel},
      series={Algorithms Combin.},
      volume={26},
      publisher={Springer, Berlin},
   },
   date={2006},
   pages={231--246},
}


\bib{mansfield1983}{article}{
   author={Mansfield, Anthony},
   title={Determining the thickness of graphs is NP-hard},
   journal={Math. Proc. Cambridge Philos. Soc.},
   volume={93},
   date={1983},
   number={1},
   pages={9--23},
}

\bib{mayhew-08}{article}{
   author={Mayhew, Dillon},
   title={Matroid complexity and nonsuccinct descriptions},
   journal={SIAM J. Discrete Math.},
   volume={22},
   date={2008},
   number={2},
   pages={455--466},
}

\bib{MNWW-11}{article}{
   author={Mayhew, Dillon},
   author={Newman, Mike},
   author={Welsh, Dominic},
   author={Whittle, Geoff},
   title={On the asymptotic proportion of connected matroids},
   journal={European J. Combin.},
   volume={32},
   date={2011},
   number={6},
   pages={882--890},
}

\bib{mayhew-welsh-13}{article}{
   author={Mayhew, Dillon},
   author={Welsh, Dominic},
   title={On the number of sparse paving matroids},
   journal={Adv. in Appl. Math.},
   volume={50},
   date={2013},
   number={1},
   pages={125--131},
}



\bib{merino-97}{article}{
   author={Merino L\'{o}pez, Criel},
   title={Chip firing and the Tutte polynomial},
   journal={Ann. Comb.},
   volume={1},
   date={1997},
   number={3},
   pages={253--259},
}


\bib{merino-welsh-99}{article}{
   author={Merino, C.},
   author={Welsh, D. J. A.},
   title={Forests, colorings and acyclic orientations of the square lattice},
   note={On combinatorics and statistical mechanics},
   journal={Ann. Comb.},
   volume={3},
   date={1999},
   number={2-4},
   pages={417--429},
}

\bib{milner-welsh-75}{article}{
   author={Milner, E. C.},
   author={Welsh, D. J. A.},
   title={On the computational complexity of graph theoretical properties},
   conference={
      title={Proceedings of the Fifth British Combinatorial Conference},
      address={Univ. Aberdeen, Aberdeen},
      date={1975},
   },
   book={
      series={Congress. Numer.},
      volume={No. XV},
      publisher={Utilitas Math., Winnipeg, MB},
   },
   date={1976},
   pages={471--487},
}


\bib{morgan-welsh1965}{article}{
   author={Morgan, R. W.},
   author={Welsh, D. J. A.},
   title={A two-dimensional Poisson growth process},
   journal={J. Roy. Statist. Soc. Ser. B},
   volume={27},
   date={1965},
   number={3},
   pages={497--504},
}
 

\bib{mosca-98}{article}{
   author={Mosca, Michele},
   title={Removing edges can increase the average number of colours in the
   colourings of a graph},
   journal={Combin. Probab. Comput.},
   volume={7},
   date={1998},
   number={2},
   pages={211--216},
}



\bib{nakasawa1935}{article}{
   author={Nakasawa, Takeo},
   title={Zur Axiomatik der linearen Abh\"angigkeit. I},
   journal={Sci. Rep. Tokyo Bunrika Daigaku Sect. A},
   volume={2},
   date={1935},
   number={43},
   pages={129--149},
}

\bib{nakasawa1936a}{article}{
   author={Nakasawa, Takeo},
   title={Zur Axiomatik der linearen Abh\"angigkeit. II},
   journal={Sci. Rep. Tokyo Bunrika Daigaku Sect. A},
   volume={3},
   date={1936},
   number={51},
   pages={17--41},
}

\bib{nakasawa1936b}{article}{
   author={Nakasawa, Takeo},
   title={Zur Axiomatik der linearen Abh\"angigkeit. III.  Schluss},
   journal={Sci. Rep. Tokyo Bunrika Daigaku Sect. A},
   volume={3},
   date={1936},
   number={55},
   pages={77--90},
}


\bib{NashWilliams1967}{article}{
   author={Nash-Williams, C. St.~J. A. },
   title={An application of matroids to graph theory},
   conference={
      title={Theory of Graphs: International Symposium (Rome, July 1966)},
      organization={International Computation Centre, Rome},
   },
   book={
      publisher={Dunod, Paris / Gordon and Breach, New York},
   },
   date={1967},
   pages={263--265},
}


\bib{noble-98}{article}{
   author={Noble, S. D.},
   title={Evaluating the Tutte polynomial for graphs of bounded tree-width},
   journal={Combin. Probab. Comput.},
   volume={7},
   date={1998},
   number={3},
   pages={307--321},
}

\bib{noble-07}{article}{
   author={Noble, Steven D.},
   title={Complexity of graph polynomials},
   conference={
      title={Combinatorics, complexity, and chance},
   },
   book={
      series={Oxford Lecture Ser. Math. Appl.},
      volume={34},
      publisher={Oxford Univ. Press, Oxford},
   },
   date={2007},
   pages={191--212},
}

\bib{noble-welsh1999}{article}{
   author={Noble, S. D.},
   author={Welsh, D. J. A.},
   title={A weighted graph polynomial from chromatic invariants of knots},
   note={Symposium \`a la M\'emoire de Fran\c{c}ois Jaeger (Grenoble, 1998)},
   journal={Ann. Inst. Fourier (Grenoble)},
   volume={49},
   date={1999},
   number={3},
   pages={1057--1087},
}


\bib{oxley-92}{book}{
   author={Oxley, James G.},
   title={Matroid theory},
   series={Oxford Science Publications},
   publisher={The Clarendon Press, Oxford University Press, New York},
   date={1992},
   pages={xii+532},
}

\bib{oxley-07}{article}{
   author={Oxley, James},
   title={The contributions of Dominic Welsh to matroid theory},
   conference={
      title={Combinatorics, complexity, and chance},
   },
   book={
      series={Oxford Lecture Ser. Math. Appl.},
      volume={34},
      publisher={Oxford Univ. Press, Oxford},
   },
   date={2007},
   pages={234--259},
}

\bib{oxley-11}{book}{
   author={Oxley, James},
   title={Matroid theory},
   series={Oxford Graduate Texts in Mathematics},
   volume={21},
   edition={2},
   publisher={Oxford University Press, Oxford},
   date={2011},
   pages={xiv+684},
}

\bib{oxley-welsh-79b}{article}{
   author={Oxley, J. G.},
   author={Welsh, D. J. A.},
   title={On some percolation results of J. M. Hammersley},
   journal={J. Appl. Probab.},
   volume={16},
   date={1979},
   number={3},
   pages={526--540},
}

\bib{oxley-welsh-79a}{article}{
   author={Oxley, J. G.},
   author={Welsh, D. J. A.},
   title={The Tutte polynomial and percolation},
   conference={
      title={Graph theory and related topics},
      address={Proc. Conf., Univ. Waterloo, Waterloo, Ont.},
      date={1977},
   },
   book={
      publisher={Academic Press, New York-London},
   },
   date={1979},
}


\bib{oxley-whittle1993a}{article}{
   author={Oxley, J. G.},
   author={Whittle, G. P.},
   title={Tutte invariants for 2-polymatroids},
   conference={
      title={Proc. AMS-IMS-SIAM Joint Summer Research Conf. on Graph Minors},
      address={University of Washington, Seattle, Washington},
      date={22 June -- 5 July 1991},
   },
   book={
      editor={Robertson, N.},
      editor={Seymour, P.},
      title={Graph structure theory},
      series={Contemp. Math.},
      volume={147},
      publisher={Amer. Math. Soc., Providence, RI},
   },
   pages={9--19},
   date={1993},
}

\bib{oxley-whittle1993b}{article}{
   author={Oxley, J. G.},
   author={Whittle, G. P.},
   title={A characterization of Tutte invariants of 2-polymatroids},
   journal={J. Combin. Theory Ser. B},
   volume={59},
   date={1993},
   number={2},
   pages={210--244},
}

\bib{pendavingh-vanderpol-15}{article}{
   author={Pendavingh, R. A.},
   author={van der Pol, J. G.},
   title={Counting matroids in minor-closed classes},
   journal={J. Combin. Theory Ser. B},
   volume={111},
   date={2015},
   pages={126--147},
}

\bib{petford-welsh-89}{article}{
   author={Petford, A. D.},
   author={Welsh, D. J. A.},
   title={A randomised $3$-colouring algorithm},
   note={Graph colouring and variations},
   journal={Discrete Math.},
   volume={74},
   date={1989},
   number={1-2},
   pages={253--261},
}

\bib{piff-welsh-71}{article}{
   author={Piff, M. J.},
   author={Welsh, D. J. A.},
   title={The number of combinatorial geometries},
   journal={Bull. London Math. Soc.},
   volume={3},
   date={1971},
   pages={55--56},
}


\bib{regan2024}{webpage}{
    author={Regan, K. W.},
    title={Dominic Welsh, 1938--2023},
    date={2024/03/01},
    url={https://rjlipton.com/2024/03/01/dominic-welsh-1938-2023/},
    note={In the blog \cite{lipton-regan}},
}


\bib{robinson-welsh-80}{article}{
   author={Robinson, G. C.},
   author={Welsh, D. J. A.},
   title={The computational complexity of matroid properties},
   journal={Math. Proc. Cambridge Philos. Soc.},
   volume={87},
   date={1980},
   number={1},
   pages={29--45},
}

\bib{rota-71}{article}{
   author={Rota, Gian-Carlo},
   title={Combinatorial theory, old and new},
   conference={
      title={Actes du Congr\`es International des Math\'{e}maticiens},
      address={Nice},
      date={1970},
   },
   book={
      publisher={Gauthier-Villars \'{E}diteur, Paris},
   },
   date={1971},
   pages={229--233},
}

\bib{rota-83}{article}{
   author={Rota, G.-C.},
   title={Book Reviews},
   journal={Advances in Mathematics},
   volume={47},
   pages={231},
   date={1983},
}

\bib{rota-harper-71}{article}{
   author={Rota, Gian-Carlo},
   author={Harper, L. H.},
   title={Matching theory, an introduction},
   conference={
      title={Advances in Probability and Related Topics, Vol. 1},
   },
   book={
      publisher={Dekker, New York},
   },
   date={1971},
   pages={169--215},
}

\bib{royle-12}{webpage}{
    author={Royle, G.},
    title={The Merino--Welsh Conjecture},
    date={2012-04-23},
    url={https://symomega.wordpress.com/2012/04/23/the-merino-welsh-conjecture/},
    accessdate={2024-04-20},
}

\bib{royle-15}{webpage}{
    author={Royle, G.},
    title={The Shameful Conjecture},
    date={2015-08-15},
    url={https://symomega.wordpress.com/2015/08/14/the-shameful-conjecture/},
    accessdate={2024-04-20},
}

\bib{seymour-welsh-75}{article}{
   author={Seymour, P. D.},
   author={Welsh, D. J. A.},
   title={Combinatorial applications of an inequality from statistical
   mechanics},
   journal={Math. Proc. Cambridge Philos. Soc.},
   volume={77},
   date={1975},
   pages={485--495},
   issn={0305-0041},
   review={\MR{0376378}},
   doi={10.1017/S0305004100051306},
}

\bib{seymour-welsh-78}{article}{
   author={Seymour, P. D.},
   author={Welsh, D. J. A.},
   title={Percolation probabilities on the square lattice},
   journal={Ann. Discrete Math.},
   volume={3},
   date={1978},
   pages={227--245},
}

\bib{sufrin2007}{webpage}{
    author={Sufrin, Bernard},
    title={Oxford University Computing Laboratory --- an informal prehistory},
    url={https://www.cs.ox.ac.uk/people/bernard.sufrin/personal/historyfortalk.pdf},
    date={June 2007},
    accessdate={2024-05-19},
}

\bib{talbot-welsh-10}{book}{
   author={Talbot, John},
   author={Welsh, Dominic},
   title={Complexity and cryptography},
   note={An introduction},
   publisher={Cambridge University Press, Cambridge},
   date={2006},
   pages={xii+292},
}

\bib{tutte1947}{article}{
   author={Tutte, W. T.},
   title={A ring in graph theory},
   journal={Proc. Camb. Phil. Soc.},
   volume={43},
   date={1947},
   pages={26--40},
}

\bib{tutte1948}{thesis}{
   author={Tutte, W. T.},
   title={An Algebraic Theory of Graphs},
   date={1948},
   organization={University of Cambridge},
   type={PhD thesis},
}

\bib{tutte1954}{article}{
   author={Tutte, W. T.},
   title={A contribution to the theory of chromatic polynomials},
   journal={Canad. J. Math.},
   volume={6},
   date={1954},
   pages={80--91},
}

\bib{vdHPvdP-22}{article}{
   author={van der Hofstad, Remco},
   author={Pendavingh, Rudi},
   author={van der Pol, Jorn},
   title={The number of partial Steiner systems and $d$-partitions},
   journal={Adv. Comb.},
   date={2022},
   pages={Paper No. 2, 23},
}


\bib{vertigan-welsh1992}{article}{
   author={Vertigan, D. L.},
   author={Welsh, D. J. A.},
   title={The computational complexity of the Tutte plane: the bipartite case},
   journal={Combin. Probab. Comput.},
   volume={1},
   date={1992},
   number={2},
   pages={181--187},
}


\bib{welsh-64}{thesis}{
   author={Welsh, D. J. A.},
   title={Topics in Stochastic Processes with special reference to First Passage Percolation Theory},
   date={1964},
   organization={University of Oxford},
   type={DPhil thesis},
}


\bib{welsh1965}{article}{
   author={Welsh, D. J. A.},
   title={Errors introduced by a PERT assumption},
   journal={Operations Res.},
   volume={13},
   date={1965},
   number={1},
   pages={141--143},
}

\bib{welsh1966}{article}{
   author={Welsh, D. J. A.},
   title={Super-critical arcs of a PERT network},
   journal={Operations Res.},
   volume={14},
   date={1966},
   number={1},
   pages={173--174},
}

\bib{welsh-69a}{article}{
   author={Welsh, D. J. A.},
   title={Euler and bipartite matroids},
   journal={J. Combinatorial Theory},
   volume={6},
   date={1969},
   pages={375--377},
}

\bib{welsh-69b}{article}{
   author={Welsh, Dominic J. A.},
   title={A bound for the number of matroids},
   journal={J. Combinatorial Theory},
   volume={6},
   date={1969},
   pages={313--316},
}

\bib{welsh-69c}{article}{
   author={Welsh, D. J. A.},
   title={Transversal theory and matroids},
   journal={Canadian J. Math.},
   volume={21},
   date={1969},
   pages={1323--1330},
}

\bib{welsh-70}{article}{
   author={Welsh, D. J. A.},
   title={Probability theory and its applications},
   conference={
      title={Mathematics Applied to Physics},
   },
   book={
      publisher={Springer-Verlag New York, Inc., New York},
   },
   date={1970},
   pages={465--561},
}

\bib{welsh-71a}{article}{
   author={Welsh, D. J. A.},
   title={Generalized versions of Hall's theorem},
   journal={J. Combinatorial Theory Ser. B},
   volume={10},
   date={1971},
   pages={95--101},
}

\bib{welsh-71b}{collection}{
   title={Combinatorial mathematics and its applications},
   booktitle={Proceedings of a Conference held at the Mathematical
   Institute, Oxford, from 7-10 July, 1969},
   editor={Welsh, D. J. A.},
   publisher={Academic Press, London-New York},
   date={1971},
   pages={x+364},
}

\bib{welsh-71c}{article}{
   author={Welsh, D. J. A.},
   title={Combinatorial problems in matroid theory},
   conference={
      title={Combinatorial Mathematics and its Applications},
      address={Proc. Conf., Oxford},
      date={1969},
   },
   book={
      publisher={Academic Press, London-New York},
   },
   date={1971},
   pages={291--306},
}

\bib{welsh-76}{book}{
   author={Welsh, D. J. A.},
   title={Matroid theory},
   series={L. M. S. Monographs},
   volume={No. 8},
   publisher={Academic Press [Harcourt Brace Jovanovich, Publishers],
   London-New York},
   date={1976},
   pages={xi+433},
   note={Dover, New York, 2010}
}

\bib{welsh-83}{article}{
   author={Welsh, D. J. A.},
   title={Randomised algorithms},
   journal={Discrete Appl. Math.},
   volume={5},
   date={1983},
   number={1},
   pages={133--145},
}

\bib{welsh-88}{book}{
   author={Welsh, Dominic},
   title={Codes and cryptography},
   series={Oxford Science Publications},
   publisher={The Clarendon Press, Oxford University Press, New York},
   date={1988},
   pages={xii+257},
}



\bib{welsh1990}{article}{
   author={Welsh, D. J. A.},
   title={The computational complexity of some classical problems from statistical physics},
   book={
      editor={Grimmett, G. R.},
      editor={Welsh, D. J. A.},
      title={Disorder in physical systems: A volume in honour of John M. Hammersley on the occasion of his 70th birthday},
      publisher={Oxford University Press},
   },
   date={1990},
   pages={307--321},
}


\bib{welsh-93a}{article}{
   author={Welsh, D. J. A.},
   title={Percolation in the random cluster process and $Q$-state Potts model},
   journal={J. Phys. A},
   volume={26},
   date={1993},
   number={11},
   pages={2471--2483},
}

\bib{welsh-93b}{book}{
   author={Welsh, D. J. A.},
   title={Complexity: knots, colourings and counting},
   series={London Mathematical Society Lecture Note Series},
   volume={186},
   publisher={Cambridge University Press, Cambridge},
   date={1993},
   pages={viii+163},
}

\bib{welsh-96}{article}{
   author={Welsh, D. J. A.},
   title={Counting colourings and flows in random graphs},
   conference={
      title={Combinatorics, Paul Erd\H{o}s is eighty, Vol. 2},
      address={Keszthely},
      date={1993},
   },
   book={
      series={Bolyai Soc. Math. Stud.},
      volume={2},
      publisher={J\'{a}nos Bolyai Math. Soc., Budapest},
   },
   date={1996},
   pages={491--505},
}

\bib{welsh1997a}{article}{
   author={Welsh, D. J. A.},
   title={Randomised approximation in the Tutte plane},
   book={
      editor={Bollob\'as, B.},
      editor={Thomason, A.},
      title={Combinatorics, geometry and probability},
      publisher={Cambridge University Press},
   },
   date={1997},
   pages={549--555},
}

\bib{welsh1997b}{article}{
   author={Welsh, D. J. A.},
   title={Approximate counting},
   conference={
      title={Proc. 16th British Combinatorial Conf.},
      address={Queen Mary and Westfield College, University of London},
      date={July 1997},
   },
   book={
      editor={Bailey, R. A.},
      title={Surveys in combinatorics, 1997},
      series={London Math. Soc. Lecture Note Ser.},
      volume={241},
      publisher={Cambridge University Press},
   },
   date={1997},
   pages={287--323},
}

\bib{welsh-98}{article}{
   author={Welsh, Dominic},
   title={Percolation and the random cluster model: combinatorial and
   algorithmic problems},
   conference={
      title={Probabilistic methods for algorithmic discrete mathematics},
   },
   book={
      series={Algorithms Combin.},
      volume={16},
      publisher={Springer, Berlin},
   },
   date={1998},
   pages={166--194},
}

\bib{welsh1999}{article}{
   author={Welsh, D. J. A.},
   title={The Tutte polynomial},
   note={Special issue: Statistical physics methods in discrete probability, combinatorics, and theoretical computer science (Princeton, NJ, 1997)},
   journal={Random Structures Algorithms},
   volume={15},
   date={1999},
   number={3--4},
   pages={210--228},
}


\bib{welsh2003}{article}{
   author={Welsh, D. J. A.},
   title={Crispin St.~J. A. Nash-Williams (1932--2001)},
   journal={Bull. London Math. Soc.},
   volume={35},
   date={2003},
   number={6},
   pages={829--844},
}

\bib{welsh2006}{webpage}{
    author={Welsh, D.},
    title={Convocation Address},
    url={https://uwaterloo.ca/combinatorics-and-optimization/news/dominic-welsh-awarded-honorary-dmath-degree-and-addressed},
    date={2006-06-16},
    accessdate={2024-05-20},
}


\bib{welsh-powell-67}{article}{
    author={Welsh, D. J. A.},
    author={Powell, M. B.},
    title = {An upper bound for the chromatic number of a graph and its application to timetabling problems},
    journal = {The Computer Journal},
    volume = {10},
    number = {1},
    pages = {85-86},
    year = {1967},
}


\bib{welsh-merino-00}{article}{
   author={Welsh, D. J. A.},
   author={Merino, C.},
   title={The Potts model and the Tutte polynomial},
   note={Probabilistic techniques in equilibrium and nonequilibrium
   statistical physics},
   journal={J. Math. Phys.},
   volume={41},
   date={2000},
   number={3},
   pages={1127--1152},
}


\bib{welsh-whittle1999}{article}{
   author={Welsh, Dominic J. A.},
   author={Whittle, Geoffrey P.},
   title={Arrangements, channel assignments, and associated polynomials},
   journal={Adv. in Appl. Math.},
   volume={23},
   date={1999},
   number={4},
   pages={375--406},
}


\bib{whitney1932}{article}{
   author={Whitney, Hassler},
   title={The coloring of graphs},
   journal={Ann. of Math. (2)},
   volume={33},
   date={1932},
   pages={688--718},
}


\bib{whitney1935}{article}{
   author={Whitney, Hassler},
   title={On the abstract properties of linear dependence},
   journal={Amer. J. Math.},
   volume={57},
   date={1935},
   pages={509--533},
}



\bib{wild-05}{article}{
   author={Wild, Marcel},
   title={The asymptotic number of binary codes and binary matroids},
   journal={SIAM J. Discrete Math.},
   volume={19},
   date={2005},
   number={3},
   pages={691--699},
}

\bib{wood-1}{webpage}{
    author={Wood, D. R.},
    title={The Academic Family Tree of Dominic Welsh (1938--2023)},
    url={https://users.monash.edu.au/~davidwo/files/Welsh-FamilyTree.pdf},
    accessdate={2024-04-20},
}

\bib{wood-2}{webpage}{
    author={Wood, D. R.},
    title={The Academic Family Tree of John M. Hammersley (1920--2004)},
    url={https://users.monash.edu.au/~davidwo/files/Hammersley-FamilyTree.pdf},
    accessdate={2024-04-20},
}

\end{biblist}

\end{bibdiv}
\end{document}